\documentclass[11pt]{article}
\usepackage{amsthm,amsmath,latexsym,amssymb}

\newcommand{\bP}{{\rm |\kern-.15em P}}
\newcommand{\Q}{\kern.3em\rule{.07em}{.65em}\kern-.3em{\rm Q}}
\newcommand{\R}{{\rm I\kern-.15em R}}
\newcommand{\D}{{\rm |\kern-.15em D}}
\newcommand{\h}{{\rm |\kern-.15em H}}
\newcommand{\C}{\kern.3em\rule{.07em}{.65em}\kern-.3em{\rm C}}
\newcommand{\T}{{\rm T\kern-.35em T}}

\theoremstyle{plain}
\newtheorem{theorem}{Theorem}[section]

\newtheorem{proposition}[theorem]{Proposition}
\newtheorem{corollary}[theorem]{Corollary}

\theoremstyle{definition}
\newtheorem{definition}[theorem]{Definition}

\theoremstyle{remark}
\newtheorem{remark}[theorem]{Remark}

\begin{document}
\title{Picard theorems for Keller mappings in dimension two and the phantom curve}
\author{Ronen Peretz}
 
\maketitle

\begin{abstract}
Let $F=(P,Q)\in\mathbb{C}[X,Y]^{2}$ be a polynomial mapping over the complex field $\mathbb{C}$. Suppose that
$$
\det\,J_{F}(X,Y):=\frac{\partial P}{\partial X}\frac{\partial Q}{\partial Y}-
\frac{\partial P}{\partial Y}\frac{\partial Q}{\partial X}=a\in\mathbb{C}^{\times}.
$$
A mapping that satisfies the assumptions above is called a Keller mapping. In this paper we estimate the size of
the co-image of $F$. We give a sufficient condition for surjectivity of Keller mappings in terms of its Phantom
curve. This curve is closely related to the asymptotic variety of $F$.
\end{abstract}

\section{Introduction}

In this paper we will prove, among other things, the following results: \\
\\
{\bf Theorem 4.1} {\it If $F\in\mathbb{C}[X,Y]^2$ satisfies $\det J_F(X,Y)\in\mathbb{C}^{\times}$, and
$\forall\,R\in R_0(F),\,\,\{ S_R(X,Y)=0\}\cap {\rm sing}(R)=\emptyset$, then $F(\mathbb{C}^2)=\mathbb{C}^2$.} \\
\\
{\bf Theorem 5.6} {\it If $F\in\mathbb{C}[X,Y]^2$ satisfies $\det J_F(X,Y)\in\mathbb{C}^{\times}$, then} 
$$
|\mathbb{C}^2-F(\mathbb{C}^2)|\le (\deg F)^3+(\deg F)^2-(\deg F).
$$
\noindent
These results are also true over certain fields $K$ different from $\mathbb{C}$.
\\
The proofs are based on a careful analysis of the asymptotic behavior of the mapping at infinity.
The set of all the asymptotic values of $F$ is called {\bf the asymptotic variety of $F$} and is denoted by
$A(F)$. If $F$ is not an automorphism of $\mathbb{C}^2$ then this is a planar algebraic curve. Otherwise
it is an empty set. Each component of $A(F)$ is a polynomial curve, i.e. it has a normal parametrization
with polynomials. Equivalently, it has a unique place on the line at infinity (in the projectivization). 
However, every such a component is not isomorphic with $\mathbb{C}$ and hence must be singular.
These are all well known results. We can refine the description of that structure. There is a finite
set of rational but not polynomial mappings which we call {\bf a geometric basis of $F$}. It is denoted
by $R_0(F)$ and contains rational mappings of the form $L\circ (X^{-\alpha},X^{\beta}Y+X^{-\alpha}\Phi(X))$,
where $L(X,Y)$ is a fixed linear invertible mapping (depending only on $F$), $\alpha\in\mathbb{Z}^+$, 
$\beta\in\mathbb{Z}^+\cup\{0\}$, $\Phi(X)\in\mathbb{C}[X]$ and $\deg\Phi<\alpha+\beta$. Moreover, the powers 
of $X$ that effectively appear in $X^{\alpha+\beta}Y+\Phi(X)$ have a gcd which equals $1$. The cardinality of 
$R_0(F)$ equals the number of components of the asymptotic variety $A(F)$. Each rational mapping $R\in R_0(F)$ 
satisfies an identity of the form $F\circ R=G_R\in\mathbb{C}[X,Y]^2$. We call that {\bf a double asymptotic identity of $F$}. 
We call the corresponding polynomial mapping $G_R$ {\bf the $R$-dual of $F$}. The irreducible component of $A(F)$ that 
corresponds to $R\in R_0(F)$ is the polynomial curve with the following normal parametrization (meaning a 
surjective parametrization), $ \{ G_R(0,Y)\,\vert\,Y\in\mathbb{C}\}$. We call this component {\bf the $R$-component of $A(F)$}. 
If its implicit representation is $H_R(X,Y)=0$ for some irreducible $H_R\in\mathbb{C}[X,Y]$ then $(H_R\circ G_R)(0,Y)\equiv 0$. 
In fact we prove that $ (H_R\circ G_R)(X,Y)=X^{\beta-\alpha}S_R(X,Y)$, where $1\le\beta-\alpha$, and where 
$S_R\in\mathbb{C}[X,Y]$. The planar algebraic curve $\{ S_R(X,Y)=0\}$, is called {\bf the $R$-phantom curve}.
If $\forall\,R\in R_0(F),\,\,\{ S_R(X,Y)=0\}\cap{\rm sing}(R)=\emptyset$ then the  mapping $F$ must be surjective! 
The reason is the following: If the finite set $\mathbb{C}^2-F(\mathbb{C}^2)\ne\emptyset$ is non-empty then any 
$(a,b)\in\mathbb{C}^2-F(\mathbb{C}^2)$ is an asymptotic value of $F$. We call such an asymptotic value, {\bf a 
Picard-exceptional value of $F$} (as is the terminology in the theory of analytic functions). The mapping $F$ is 
surjective if and only if it has no Picard-exceptional values. Let $R$ be an element of $R_0(F)$ that corresponds to $(a,b)$. 
This means that $(a,b)\in \{H_R(X,Y)=0\}-F(\mathbb{C}^2)$. As explained above we have 
$\{H_R(X,Y)=0\}=\{G_R(0,Y)\,|\,Y\in\mathbb{C}\}$. The inverse image $G_R^{-1}(\{H_R(X,Y)=0\}$
equals the union $\{X=0\}\cup\{S_R(X,Y)=0\}$, i.e. the union of the singular set of $R$ and the $R$-phantom curve.
If these two sets are disjoint then $R$ is defined on every point of the $R$-phantom curve $\{S_R(X,Y)=0\}$. In particular
$G_R(\{S_R(X,Y)=0\})=F(R(\{S_R(X,Y)=0\}))\subseteq F(\mathbb{C}^2)$. This means that if $F$ has Picard-exceptional
values on $\{H_R(X,Y)=0\}$ they must belong to the difference set $\{H_R(X,Y)=0\}-G_R(\{S_R(X,Y)=0\})$.
We will prove that $\{H_R(X,Y)=0\}-G_R(\{S_R(X,Y)=0\})\subseteq \{F(0,Y)\,|\,Y\in\mathbb{C}\}$ and so the $R$-component
$\{H_R(X,Y)=0\}$ of the asymptotic variety $A(F)$ contains no Picard-exceptional values of $F$. Since this is true
for any $R\in R_0(F)$, it follows that $F$ has no Picard-exceptional values. Hence $F$ is a surjective mapping.
We will prove that for certain choices of the parameters $\alpha$ and $\beta$ the condition 
$\{ S_R(X,Y)=0\}\cap{\rm sing}(R)=\emptyset$ is satisfied. As a corollary we are able to prove that: \\
{\it If $N\in\mathbb{Z}^+$ and $a_1,\ldots,a_N\in\mathbb{C}$ then $I((X^{-N},X^{N+1}Y+a_N X^N+\ldots a_1 X))$ contains
no Jacobian pair.} \\
This result should be compared to the main theorems in \cite{rp,rp1,w} which handled other families of subalgebras
of $\mathbb{C}[X,Y]$. The case $N=1$ is in the intersection of the two families and originally was proved
by L. Makar-Limanov using techniques from weighted graded algebras. We recall that Pinchuk's counterexample to 
the Real Jacobian Conjecture is contained in the real version $\mathbb{R}[V,VU,VU^2+U]$ of the case $N=1$.

\section{The structure of the asymptotic variety}

{\bf An $F$-asymptotic value} $(a,b)\in\mathbb{C}^2$ is a limiting value of $F:\,\mathbb{C}^2\rightarrow\mathbb{C}^2$
along a smooth curve that tends to infinity. The set of all the $F$-asymptotic values is called 
{\bf the $F$-asymptotic variety} and it is denoted by $A(F)$. Any smooth curve as above is called {\bf an
asymptotic tract of $F$ that corresponds to $(a,b)$}.

\begin{theorem}
If $F$ is a two dimensional Keller mapping, and if $F$ satisfies {\bf the $Y$-degree condition}:
$\deg F=\deg_Y P=\deg_Y Q$, $(F=(P,Q)$). Then $\forall\,(a,b)\in A(F)$ there exists a rational
mapping: $R(X,Y)=(X^{-\alpha},X^{\beta}Y+X^{-\alpha}\phi(X))$, ($\alpha\in\mathbb{Z}^+,\,\beta\in
\mathbb{Z}^+\cup\{0\},\,\phi(X)\in\mathbb{C}[X],\,\deg\phi<\alpha+\beta$ and the gcd of the powers of $X$ that
effectively appear in $X^{\alpha+\beta}Y+\phi(X)$ is $1$) with the following two properties: \\
(i) $\exists\,G_R\in\mathbb{C}[X,Y]^2$ such that $G_R=F\circ R$ off ${\rm sing}(R)=\{X=0\}$. (The
mapping $G_R$ is called {\bf the $R$-dual of $F$}) \\
(ii) $\exists\,Y_0\in\mathbb{C}$ such that $(a,b)=G_R(0,Y_0)$. \\
Moreover, the set of all those rational mappings $R$ can be chosen to be a finite set. 
\end{theorem}
\noindent
{\bf Proof.} \\
We consider the transformed  mapping:
$$
F_1(U,V)=F\left(\frac{1}{U},\frac{V}{U}\right)=\frac{A(U,V)}{U^{N}},
$$
where $N=\deg F$, $A(U,V)\in\mathbb{C}[U,V]^2$, $\deg A(U,V)=\deg A(0,V)=N$. Any asymptotic value of $F$ 
is a limiting value of $F_1(U,V)$ as $U\rightarrow 0$ and $V$ remains bounded. 
We use the following representation:
$$
F_1(U,V)=\frac{\sum_{j=0}^{N} A_{j}(V)U^{j}}{U^{N}}.  \eqno(*)
$$
In this notation $A_j(V)\in\mathbb{C}[V]^2$, $0\le j\le N$. These are the two polynomial coefficients of $U^j$.
When $U\rightarrow 0$ and $V$ remains bounded we will have $F_1(U,V)\rightarrow
\infty$ unless $V$ tends to a zero of the free term $A_{0}(V)$. This means a simultaneous zero of the 
pair of polynomials $A_0(V)$. The reason is that if $V$ is bounded away from zeros of $A_{0}(V)$, 
then the term $A_{0}(V)/U^{N}$ dominates the finite sum $F_1(U,V)$ (in equation (*)) when $U\rightarrow 0$.
We conclude that in order to determine all the possible
asymptotic values of $F$, we should consider the limits of $F_1(U,V)$ when
$U\rightarrow 0$ and $V\rightarrow a_0$, where $a_0$ is a zero of $A_{0}$. 
Let us represent the coefficients of $A$ as follows, $ A_{j}(V)=(V-a_0)^{p_{j}}B_{j}(V),\,\,\,\,\,\,0\le j\le N$.
Where if $A_{j}\equiv (0,0)$ we agree that $p_{j}=\infty $ and otherwise $p_{j}\ge 0$
and $B_{j}(a_0)\ne (0,0)$. Next we make the transformation, $V=a_0+WU^{p}$,
where $p$ will be a positive number that will be determined soon (in fact it is going
to be rational). Putting together everything we obtain,
$$
F_1(U,V)=\frac{\sum_{j=0}^{N} B_{j}(V)W^{p_{j}}U^{pp_{j}+j}}{U^{N}}.
$$
We want to determine the desired value of $p$ that will lead to a finite limiting value
of $F_1(U,V)$. If $p$ is a small positive number, then $\forall\,1\le j\le N$, $pp_{0}<pp_{j}+j$.
Thus the term that will dominate the finite sum that represents $F_1(U,V)$ is $B_{0}(a_0)W^{p_{0}}/
U^{N-pp_{0}}$. If $pp_{0}=N$ and $pp_{0}<pp_{j}+j$ for $1\le j\le N$, then
with this choice of $p=N/p_{0}$ we obtain the family of limiting values
$B_{0}(a_0)W^{p_{0}}$. In the complementary case, we choose $p>0$ so 
that $pp_{0}\le pp_{j}+j$ for $1\le j\le N$ with equality for at least
one value $j_{0}$ of $j$. We can write this choice of a value for $p$ as follows,
$$
p=\min\left\{\frac{j}{p_{0}-p_{j}}\,|\,1\le j\le N,\,p_{0}>p_{j}\right\}.
$$
We substitute this and obtain:
$$
F_1(U,a_0+WU^{p})=\frac{C(U,W)}{U^{N-pp_{0}}},
$$
where $C(U,W)$ is a mapping with coordinates that are polynomial in $U$ and $U^{p}$. It is also a polynomial in
$W$ of degree $N$ or less. Finally $C(0,W)$ is a polynomial in $W$ of degree 
$p_{0}$, which contains only the powers $W^{p_{j}}$ that satisfy the condition
$pp_{0}=pp_{j}+j$. Let $p=b_0/c$ where $b_0$ and $c$ are co prime positive integers.
Then to avoid fractional powers we make further the substitution $U=Z^{c}$.
We get,
$$
F_1(Z^{c},a_0+WZ^{b_0})=\frac{D(Z,W)}{Z^{cN-b_0p_{0}}}.   \eqno(**)
$$
Here $D(Z,W)$ is a polynomial pair in $Z$ and $W$ and $D(0,W)=C(0,W)$ is a polynomial
pair in $W$ of degree $p_{0}$. We noted above that the powers $W^{p_{j}}$ that
appear in $D(0,W)$ satisfy $pp_{0}=pp_{j}+j$. Hence for these $j$'s
$$
p_{0}-p_{j}=\frac{c\cdot j}{b_0}.
$$
Since $b_0$ and $c$ are co prime and since $p_{0}-p_{j}\in\mathbb{Z}$, $b_0$ must be a
divisor of $j$, so that the difference $p_{0}-p_{j}$ is a multiple of $c$. We
conclude that $D(0,W)$ is of the form $D(0,W)=W^{\alpha}E(W^{c})$, where $\alpha\in\mathbb{Z}^{+}\cup\{ 0\}$, 
and $E(X)$ is a polynomial pair in $X$ of degree 1 or more. \\
We see that the form of $F_1$ after our transformations, as given in equation
(**) is of the same type in the indeterminates $Z$ and $W$ as it was in equation
(*) in the indeterminates $U$ and $V$. Thus we can repeat the sequence of
transformations with suitable parameters. We conclude that any limiting value
of $F_1$ will be attained along a curve of the form $V=a_0+WU^{p}$, where 
$U\rightarrow 0$ and $W$ being bounded and $p$ is a positive rational number.
If the parameter $p$ is chosen to be smaller than our minimum value formula
then $F_1\rightarrow\infty $ along the curve. So finite limiting values for
$F_1$ can be achieved only for this value or larger ones for $p$. So with that
minimum chosen value for $p$, asymptotic values of $F_1$ that correspond to
the zero $a_0$ of the pair $A_{0}$ are achieved when $Z\rightarrow 0$ and $W$ being bounded.
If the value of $p$ is chosen to be larger than the minimum then, in fact
$W\rightarrow 0$. On repetition of the process, each asymptotic value will 
correspond to a zero $a_{1}$ of $D(0,W)$. If the multiplicity of the zero $a_{1}$ 
is $q_{0}$, then the transformation we apply is $W=a_{1}+SZ^{q}$, where $q$
is chosen using our minimum process (as with $p$). We have $q_{0}\le p_{0}
\le N$. The case $q_{0}=p_{0}$ implies $D(0,W)=d\cdot (W-a_{1})^{p_{0}}$ with a non
zero constant $d$. So $D(0,W)$ contains all the powers of $W$ from 0
to $p_{0}$. Hence the denominator of $q$, $c=1$ and the power of $Z$ in
the denominator is $N-p_{0}$ or less. We conclude, that repeating the process
will in every iteration either strictly reduce the degree of the leading
term, or strictly reduce the power of the denominator by a positive integer.
We conclude that the process must terminate. The way it will terminate is 
as follows. If we use the notation $A(U,V)/U^{N}$, the final transformation 
must be of the form $V=a_0+WU^{p}$ with $p=N/p_{0}$. The other possibility is that 
this choice of value for $p$ coincides with the minimum $\min \{j/(p_{0}-p_{j})\}$.
We obtain a curve of asymptotic values $D(0,W)$, a polynomial in $W$ of
degree $p_{0}\ge 1$. For as the denominator in the last iteration becomes 1,
our sequence of transformations assigns a new polynomial $D(Z,W)$ to the
original $F(X,Y)$. Clearly if we perform the process to every zero of the leading
term at every stage we will obtain all the asymptotic values of $F$.
Clearly we obtain in this way finitely many rational mappings of the type prescribed in
the statement of our theorem. We obtain those mappings by composing the sequence of the transformations (and invert
the result). In general asymptotic values of $F$ will be achieved as limits of $F$
along more than one of the rational curves in our construction. $\qed $ \\

\begin{remark}
We did not use the Jacobian condition in our proof. This theorem is valid for polynomial mappings
that are not necessarily Keller.
\end{remark}

\begin{theorem}
Let $F\in\mathbb{C}[X,Y]^2$ be a Keller mapping which is not an automorphism. Then there exists a 
non empty set, $R_0(F)$ of rational mappings $R\in\mathbb{C}(X,Y)^2-\mathbb{C}[X,Y]^2$ that satisfies: \\
(i) The cardinality of $R_0(F)$ equals the number of the irreducible components of the $F$-asymptotic variety $A(F)$. \\
(ii) $\forall\,R\in R_0(F)$, the set $\{G_R(a,b)\,|\,(a,b)\in {\rm sing}(R)\}$ is an irreducible component of $A(F)$. \\
(iii) We have the representation $A(F)=\bigcup_{R\in R_0(F)} \{G_R(a,b)\,|\,(a,b)\in {\rm sing}(R)\}$. \\
(iv) $\forall\,R\in R_0(F)$, the mapping $G_R=F\circ R$ has a polynomial extension to $\mathbb{C}^2$.
\end{theorem}
\noindent
{\bf Proof.} \\
We choose a linear invertible mapping $L(U,V)$ so that $F\circ L$ satisfies the $V$-degree condition.
Then the construction outlined in the proof of Theorem 2.1 gives us finitely many rational non polynomial mappings
of the form $R(X,Y)=L(X^{-\alpha},X^{\beta}Y+X^{-\alpha}\phi(X))$. This aggregate of mappings satisfy $G_R=F\circ R\in
\mathbb{C}[X,Y]^2$ off $\{X=0\}$. This proves part (iv). Any asymptotic value  of $F$ equals $G_R(0,Y)$ for some
$R$ and $Y\in\mathbb{C}$. Thus proving part (iii). The polynomial parametrization $\{G_R(0,Y)\,|\,Y\in\mathbb{C}\}$ is 
a normal parametrization (i.e. a surjective one) of one of the irreducible components of $A(F)$ and we conclude the proof 
of Theorem 2.3. $\qed $ \\

\begin{definition}
If $F\in\mathbb{C}[X,Y]^2$, then any set $R_0(F)$ that satisfies properties (i)-(iv) of Theorem 2.3 is called 
{\bf a geometric basis of $F$}.
\end{definition}

\begin{theorem}
If $F\in\mathbb{C}[X,Y]^2$ is a Keller mapping which is not an automorphism, then, \\
(i) Any irreducible component of $A(F)$ is a polynomial curve. \\
(ii) It is normally parametrized by $\{ G_R(0,Y)\,|\,Y\in\mathbb{C}\}$, where $G_R$ is
the $R$-dual of $F$ for some $R\in R_0(F)$. \\
(iii) It is a singular plane algebraic curve. \\
(iv) If $H_R(X,Y)=0$ is an implicit representation of that irreducible component of $A(F)$
($H_R\in\mathbb{C}[X,Y]$ is irreducible), then $(H_R\circ G_R)(X,Y)=X^{\gamma(R)}S_R(X,Y)$,
where $\gamma(R)\in\mathbb{Z}^+$ satisfies the inequality $\gamma(R)\le\beta-\alpha$ 
($R(X,Y)=L(X^{-\alpha},X^{\beta}Y+X^{-\alpha}\phi(X))$ as in Theorem 2.3). Also $S_R(X,Y)\in\mathbb{C}[X,Y]$.
\end{theorem}
\noindent
{\bf Proof.} \\
Parts (i) and (ii) follow by Theorem 2.3 part (ii) (a polynomial curve is a curve that has a normal 
polynomial parametrization). Now for the proof of part (iii): As implied by parts (i)-(ii) each
irreducible component of $A(F)$ is a surjective polynomial image of $\mathbb{C}$, ($\{ X=0\}$).
On the other hand by a result of N. Van Chau, Theorem 4.4 in \cite{ch1} (also a refined version in \cite{ch3}),
each such a component can not be isomorphic to $\mathbb{C}$. Since the only non singular irreducible
plane algebraic curves are isomorphic images of $\mathbb{C}$ it follows that each such a
component of $A(F)$ must be a singular plane algebraic curve. The proof of part (iv): By part (ii) we 
have $(H_R\circ G_R)(0,Y)\equiv 0$ and so $\exists\,S_R(X,Y)\in\mathbb{C}[X,Y]$ and a $\gamma(R)\in\mathbb{Z}^+$
such that $(H_R\circ G_R)(X,Y)=X^{\gamma(R)}S_R(X,Y)$ and $S_R(X,Y)\not\equiv 0$. At this point it is convenient
to note that $\beta-\alpha-1\ge 0$ and that $\det J_{G_R}(X,Y)\equiv c\cdot X^{\beta-\alpha-1}$ for
some $c\in\mathbb{C}^{\times}$. This follows by $G_R=F\circ R$ off ${\rm sing}(R)$ and by the 
Jacobian condition (satisfied by $F$). We denote $G_R=(G_1,G_2)\in\mathbb{C}[X,Y]^2$. Using the identity
$(H_R\circ G_R)(X,Y)=X^{\gamma(R)}S_R(X,Y)$ and the above note we deduce that:
$$
\left(\frac{\partial H_R}{\partial X}\circ G_R\right)(X,Y)=\left(-\frac{1}{\alpha}\right)X^{\alpha+\gamma(R)-\beta+1}
\left(\frac{\partial S_R}{\partial X}\frac{\partial G_2}{\partial Y}-\frac{\partial S_R}{\partial Y}
\frac{\partial G_2}{\partial X}\right)-
$$
$$
-\left(\frac{\gamma(R)}{\alpha}\right)X^{\alpha+\gamma(R)-\beta}\frac{\partial G_2}{\partial Y}\cdot S_R(X,Y),
$$
$$
\left(\frac{\partial H_R}{\partial Y}\circ G_R\right)(X,Y)=\left(-\frac{1}{\alpha}\right)X^{\alpha+\gamma(R)-\beta+1}
\left(\frac{\partial G_1}{\partial X}\frac{\partial S_R}{\partial Y}-\frac{\partial G_1}{\partial Y}
\frac{\partial S_R}{\partial X}\right)+
$$
$$
+\left(\frac{\gamma(R)}{\alpha}\right)X^{\alpha+\gamma(R)-\beta}\frac{\partial G_1}{\partial Y}\cdot S_R(X,Y).
$$
If $\alpha+\gamma(R)-\beta>0$ then we obtain,
$$
\forall\,Y\in\mathbb{C},\,\,(H_R\circ G_R)(0,Y)=
\left(\frac{\partial H_R}{\partial X}\circ G_R\right)(0,Y)=\left(\frac{\partial H_R}{\partial Y}\circ G_R\right)(0,Y)=0,
$$
which is not possible for an irreducible curve such as $H_R(X,Y)=0$ (with normal parametrization
$\{ G_R(0,Y)\,|\,Y\in\mathbb{C}\}$). We conclude that $\alpha+\gamma(R)-\beta\le 0$. $\qed $ \\
\\
To clarify further the relations among the parameters $\alpha,\,\beta$ and $\gamma(R)$ (after part (iv) of Theorem 2.5) 
we add the following proposition. 
\begin{proposition}
More relations among the parameters $\alpha,\,\beta$ and $\gamma(R)$ are given by: \\
(i) If $\gamma=1$, then $\beta-\alpha-1=0$ and ${\rm sing}(H_R(X,Y)=0)=\{ S_R(X,Y)=0\}\cap {\rm sing}(R)$. \\
(ii) If $\gamma\ge 2$, then $\beta-\alpha-1>0$ \\
In particular we deduce that if $\{ S_R(X,Y)=0\}\cap {\rm sing}(R)=\emptyset$ then $\gamma\ge 2$ and $\beta-\alpha-1>0$.
\end{proposition}
\noindent
{\bf Proof.} \\
We start with the identity $(H_R\circ G_R)(X,Y)=X^{\gamma(R)}S_R(X,Y)$. We use the notation $G_R=(G_1,G_2)$ and 
$H_R(U,V)$ and obtain:
$$
\left\{\begin{array}{l}
\frac{\partial H_R}{\partial U}(G_R)\cdot\frac{\partial G_1}{\partial X}+\frac{\partial H_R}{\partial V}(G_R)\cdot
\frac{\partial G_2}{\partial X}=\gamma(R)X^{\gamma(R)-1}S_R(X,Y)+X^{\gamma(R)}\frac{\partial S_R}{\partial X}, \\ \\
\frac{\partial H_R}{\partial U}(G_R)\cdot\frac{\partial G_1}{\partial Y}+\frac{\partial H_R}{\partial V}(G_R)\cdot
\frac{\partial G_2}{\partial Y}=X^{\gamma(R)}\frac{\partial S_R}{\partial Y}\,\,\,\,. \end{array}\right.
$$
To prove part (i) we use the assumption $\gamma(R)=1$ and we substitute in the system above
$X=0$. The resulting system is:
$$
\left\{\begin{array}{l}
\frac{\partial H_R}{\partial U}(G_R)\cdot\frac{\partial G_1}{\partial X}(0,Y)+\frac{\partial H_R}{\partial V}(G_R)\cdot
\frac{\partial G_2}{\partial X}(0,Y)=S_R(0,Y)\not\equiv 0, \\ \\
\frac{\partial H_R}{\partial U}(G_R)\cdot\frac{\partial G_1}{\partial Y}(0,Y)+\frac{\partial H_R}{\partial V}(G_R)\cdot
\frac{\partial G_2}{\partial Y}(0,Y)=0\,\,\,\,. \end{array}\right.
$$
Recalling the note we had in the proof of part (iii) of Theorem 2.5: $\det J_{G_R}(X,Y)=-\alpha X^{\beta-\alpha-1}$ and so 
depending whether $\beta-\alpha-1=0$ or $>0$ we conclude that the matrix
$$
\left(\begin{array}{cc} \frac{\partial G_1}{\partial X}(0,Y) & \frac{\partial G_2}{\partial X}(0,Y) \\
\frac{\partial G_1}{\partial Y}(0,Y) & \frac{\partial G_2}{\partial Y}(0,Y)\end{array}\right),
$$
is either always invertible (when $\beta-\alpha-1=0$) or always not invertible (when $\beta-\alpha-1>0$). In our case
the matrix $J_{G_R}(0,Y)^T$ can not always be not invertible because this would mean that the two equations in the second
system above are proportional which is an absurd because $S_R(0,Y)\not\equiv 0$. Thus in case $\gamma=1$ we must have
$\beta-\alpha-1=0$. So the second system has a unique solution $\forall\,Y\in\mathbb{C}$ and in particular,
$$
\frac{\partial H_R}{\partial U}(G_R(0,Y))=\frac{\partial H_R}{\partial V}(G_R(0,Y))=0\Leftrightarrow S_R(0,Y)=0\Leftrightarrow
$$
$$
\Leftrightarrow (0,Y)\in\{S_R(X,Y)=0\}\cap {\rm sing}(R).
$$
To prove part (ii) we use the assumption $\gamma(R)\ge 2$ and we substitute in the first system $X=0$.
We obtain the following system:
$$
\left\{\begin{array}{l}
\frac{\partial H_R}{\partial U}(G_R)\cdot\frac{\partial G_1}{\partial X}(0,Y)+\frac{\partial H_R}{\partial V}(G_R)\cdot
\frac{\partial G_2}{\partial X}(0,Y)=0, \\ \\
\frac{\partial H_R}{\partial U}(G_R)\cdot\frac{\partial G_1}{\partial Y}(0,Y)+\frac{\partial H_R}{\partial V}(G_R)\cdot
\frac{\partial G_2}{\partial Y}(0,Y)=0\,\,\,\,. \end{array}\right.
$$
In this case the matrix $J_{G_R}(0,Y)^T$ can not always be invertible, because it would imply that
$$
H_R(G_R(0,Y))\equiv\frac{\partial H_R}{\partial U}(G_R(0,Y))\equiv\frac{\partial H_R}{\partial V}(G_R(0,Y))\equiv 0.
$$
Hence, in this case we must have $\beta-\alpha-1>0$. In particular if $\{ S_R(X,Y)=0\}\cap {\rm sing}(R)=\emptyset$ and
$\gamma=1$, then by (i) ${\rm sing}(H_R(X,Y)=0)=\emptyset$ which is a contradiction. Hence, $\gamma\ge 2$. $ \qed $ \\

\begin{remark}
On the next section we will prove a more accurate version of Theorem 2.5 (iv).
\end{remark}

\section{The relation $\gamma=\beta-\alpha$, and the geometry of the $R$-phantom curve}

\begin{theorem}
Let $F$ be a Keller mapping which is not a $\mathbb{C}^2$ automorphism. Assume that $F$ satisfies
the $Y$-degree condition. Then $\forall\,R(X,Y)=(X^{-\alpha},X^{\beta}Y+X^{-\alpha}\phi(X))\in R_0(F)$ we
have the identity $H_R(G_R(X,Y))=X^{\beta-\alpha}S_R(X,Y)$.
\end{theorem}
\noindent
{\bf Proof.} \\
Let $G_R=(G_1,G_2)=F\circ R$ (off ${\rm sing}(R)$) be the $R$-dual of $F$. Let $H_R(X,Y)=0$ be the $R$-component of $A(F)$. Then
$H_R(G_R(X,Y))=X^{\gamma}S_R(X,Y)$ for some $\gamma\in\mathbb{Z}^+$, $S_R(X,Y)\in\mathbb{C}[X,Y]$, $S_R(0,Y)\not\equiv 0$.
This follows by Hilbert's Nullstellensatz and the irreducibility of $H_R$ (where $\gamma$ absorbs all the $X$-powers). So
$$
\left\{\begin{array}{l}
\frac{\partial H_R}{\partial U}(G_R)\cdot\frac{\partial G_1}{\partial X}+\frac{\partial H_R}{\partial V}(G_R)\cdot
\frac{\partial G_2}{\partial X}=\gamma X^{\gamma-1}S_R(X,Y)+X^{\gamma}\frac{\partial S_R}{\partial X}, \\ \\
\frac{\partial H_R}{\partial U}(G_R)\cdot\frac{\partial G_1}{\partial Y}+\frac{\partial H_R}{\partial V}(G_R)\cdot
\frac{\partial G_2}{\partial Y}=X^{\gamma}\frac{\partial S_R}{\partial Y}\,\,\,\,. \end{array}\right.
$$
We think of this system as a $2\times 2$ linear system in the two unknowns $(\partial H_R/\partial U)(G_R)$ and
$(\partial H_R/\partial V)(G_R)$. The coefficients matrix is
$$
\left(\begin{array}{cc} \frac{\partial G_1}{\partial X} & \frac{\partial G_2}{\partial X} \\
\frac{\partial G_1}{\partial Y} & \frac{\partial G_2}{\partial Y}\end{array}\right)=J_{G_R}(X,Y)^T.
$$
The determinant of this matrix is $\det J_{G_R}(X,Y)^T=-\alpha X^{\beta-\alpha-1}$. By Cramer's Rule we have,
$$
-\alpha X^{\beta-\alpha-\gamma}\frac{\partial H_R}{\partial U}(G_R)=
\left|\begin{array}{cc} \gamma S_R+X(\partial S_R/\partial X) & \partial G_2/\partial X \\
X(\partial S_R/\partial Y) & \partial G_2/\partial Y \end{array}\right|,
$$
and 
$$
-\alpha X^{\beta-\alpha-\gamma}\frac{\partial H_R}{\partial V}(G_R)=
\left|\begin{array}{cc} \partial G_1/\partial X & \gamma S_R+X(\partial S_R/\partial X) \\
\partial G_1/\partial Y & X(\partial S_R/\partial Y) \end{array}\right|.
$$
We evaluate for $X=0$:
$$
-\alpha 0^{\beta-\alpha-\gamma}\frac{\partial H_R}{\partial U}(G_R(0,Y))=\gamma S_R(0,Y)\frac{\partial G_2}{\partial Y}(0,Y),
$$
$$
-\alpha 0^{\beta-\alpha-\gamma}\frac{\partial H_R}{\partial V}(G_R(0,Y))=-\gamma S_R(0,Y)\frac{\partial G_1}{\partial Y}(0,Y).
$$
We consider the second equation and recall that the specialization $X=0$ is an operator that commutes
with $\partial/\partial Y$. Hence $(\partial G_1/\partial Y)(0,Y)=d G_1(0,Y)/dY$ and so if $(\partial G_1/\partial Y)(0,Y)
\equiv 0$ then $G_1(0,Y)\equiv c$ a constant. So in this case the curve $\{(G_1(0,Y),G_2(0,Y))\,|\,Y\in\mathbb{C}\}$
is either a point or $\{(c,Y)\,|\,Y\in\mathbb{C}\}$. But both possibilities can not occur because
this curve is the $R$-component of $A(F)$, $H_R(X,Y)=0$ which is a singular non-degenerate planar algebraic curve.
We conclude that $(\partial G_1/\partial Y)(0,Y)\not\equiv 0$ and so $-\gamma S_R(0,Y)(\partial G_1/\partial Y)(0,Y)\not\equiv 0$.
Thus $-\alpha 0^{\beta-\alpha-\gamma}\frac{\partial H_R}{\partial V}(G_R(0,Y))\not\equiv 0$ which proves
that $\beta-\alpha-\gamma=0$. $\qed $ \\

\begin{corollary}
Every intersection point of ${\rm sing}(R)$ and the $R$-phantom $S_R(X,Y)=0$ has a $G_R$-image which is a singular point
of the $R$-component of $A(F)$, $H_R(X,Y)=0$.
\end{corollary}
\noindent
{\bf Proof.} \\
Using the proof of Theorem 3.1 we get:
$$
-\alpha \frac{\partial H_R}{\partial U}(G_R(0,Y))=(\beta-\alpha) S_R(0,Y)\frac{\partial G_2}{\partial Y}(0,Y),
$$
$$
-\alpha \frac{\partial H_R}{\partial V}(G_R(0,Y))=-(\beta-\alpha) S_R(0,Y)\frac{\partial G_1}{\partial Y}(0,Y).
$$
This proves that 
$$
S_R(0,Y)=0\Longrightarrow \frac{\partial H_R}{\partial U}(G_R(0,Y))=\frac{\partial H_R}{\partial V}(G_R(0,Y))=0.\,\,\,\qed
$$
\begin{corollary}
$G_R({\rm sing}(S_R=0))\cup G_R(\{S_R=0\}\cap {\rm sing}(R))={\rm sing}(H_R=0)$.
\end{corollary}
\noindent
{\bf Proof.} \\
If $(X_0,Y_0)$ is a singular point of the $R$-phantom curve which is off ${\rm sing}(R)$, then $X_0\ne 0$,
$S_R(X_0,Y_0)=(\partial S_R/\partial X)(X_0,Y_0)=(\partial S_R/\partial Y)(X_0,Y_0)=0$. Hence also
$(\partial H_R/\partial U)(G_R(X_0,Y_0))=(\partial H_R/\partial V)(G_R(X_0,Y_0))=0$. This follow by the determinantial
formulas in the proof of Theorem 3.1. Hence $G_R(X_0,Y_0)$ is a singular point of $H_R(X,Y)=0$. This and Corollary 3.2
prove that $G_R({\rm sing}(S_R=0))\cup G_R(\{S_R=0\}\cap {\rm sing}(R))\subseteq{\rm sing}(H_R=0)$. The singular locus 
of ${\rm sing}(H_R=0)$ contains also points $(G_1(0,Y_0),G_2(0,Y_0))$ for which 
$(\partial G_1/\partial Y)(0,Y_0)=(\partial G_2/\partial Y)(0,Y_0)=0$. If such a singular point $(G_1(0,Y_0),G_2(0,Y_0))$ 
coincides with $(G_1(a,b),G_2(a,b))$ for which $S_R(a,b)=0$ then if $a\ne 0$ $G_R$ is a local diffeomorphism at $(a,b)$ 
which implies that $(a,b)$ is also a singular point of $S_R(X,Y)=0$. But such a point was already counted for on the 
set on the left hand side. $\qed $ \\
\\ 
A very important fact that we would like to point out in the result on the next section is that the disjointness
of the singular locus of $R$ and the $R$-phantom curve implies the surjectivity of the mapping $F$. The next
result proves that this disjointness holds true in the special case $\beta=\alpha+1$,

\begin{theorem}
Let $F$ be a Keller mapping which is not a $\mathbb{C}^2$ automorphism. Then $\forall\,R\in R_0(F)$  of the form
$R(X,Y)=L\circ (X^{-\alpha},X^{\alpha+1}Y+X^{-\alpha}\phi(X))$ we have ${\rm sing}(R)\cap\{ S_R(X,Y)=0\}=\emptyset$.
\end{theorem}
\noindent
{\bf Proof.} \\
Without losing the generality we may assume that $F$ satisfies the $Y$-degree condition. This implies that
each $R\in R_0(F)$ could be chosen to have the following form: $R(X,Y)=(X^{-\alpha},X^{\beta}Y+X^{-\alpha}\phi(X))$,
$\alpha,\beta\in\mathbb{Z}^+$, $\alpha<\beta$, $\phi(X)\in\mathbb{C}[X]$, $\deg\phi<\alpha+\beta$ and the gcd
of the set of $X$-powers that effectively appear in $X^{\alpha+\beta}Y+\phi(X)$ equals to $1$. 
We assume that $\beta=\alpha+1$. In this case we have $\det J_R(X,Y)=-\alpha,\,\,\forall\,(X,Y)\not\in {\rm sing}(R)$.
By the relation $F\circ R=G_R$, the $R$-dual of $F$, it follows that $\det J_{G_R}\in\mathbb{C}^{\times}$ (since
$F$ is Keller). Thus in this case $G_R$ is Keller as well. We know that the pre-image  of the $R$-irreducible component
of $A(F)$ by $G_R$ equals the union $\{S_R(X,Y)=0\}\cup {\rm sing}(R)$. In other words we have 
$G_R^{-1}(\{H_R(X,Y)=0\})=G_R^{-1}(\{G_R(0,Y)\,|\,Y\in\mathbb{C}\})=\{S_R(X,Y)=0\}\cup {\rm sing}(R)$. Now let us assume, 
in order to get a contradiction that  ${\rm sing}(R)\cap\{ S_R(X,Y)=0\}\ne\emptyset$. Say 
$(a,b)\in {\rm sing}(R)\cap\{ S_R(X,Y)=0\},\,\,(a=0)$. Then there exist two sequences $(a_n^1,b_n^1)\in{\rm sing}(R)$, 
$(a_n^2,b_n^2)\in\{ S_R(X,Y)=0\}$ so that: \\
(1) $\lim (a_n^1,b_n^1)=\lim (a_n^2,b_n^2)=(a,b)$. \\
(2) $\forall\,n,\,\,(a_n^1,b_n^1)\ne (a_n^2,b_n^2)$. \\
(3) $\forall\,n,\,\, G_R(a_n^1,b_n^1)=G_R(a_n^2,b_n^2)$. \\
Hence $(a,b)$ is singular point of the mapping $G_R(X,Y)$ and in particular $\det J_{G_R}(a,b)=0$. This contradicts
the fact that in our case $\beta=\alpha+1$, and as explained above this implies that $G_R(X,Y)$ is Keller, i.e.
$\det J_{G_R}\in\mathbb{C}^{\times}$. This completes the proof of the theorem. $ \qed $ \\

\begin{corollary}
If $N\in\mathbb{Z}^+$ and $a_1,\ldots,a_N\in\mathbb{C}$ then $I((X^{-N},X^{N+1}Y+a_N X^N+\ldots a_1 X))$ contains
no Jacobian pair.
\end{corollary}
\noindent
{\bf Proof.} \\
Suppose to the contrary that $F\in I((X^{-N},X^{N+1}Y+a_N X^N+\ldots a_1 X))$ is a Keller mapping. Then 
$R(X,Y)=(X^{-N},X^{N+1}Y+a_N X^N+\ldots a_1 X)\in R_0(F)$ and by Proposition 2.6(i) we have 
${\rm sing}(H_R(X,Y)=0)=\{ S_R(X,Y)=0\}\cap {\rm sing}(R)$. On the other hand, by case 2 in the proof of
Theorem 3.4, this implies that ${\rm sing}(H_R(X,Y)=0)=\emptyset$. This contradicts the fact that 
the $R$-component of $A(F)$, $\{ H_R(X,Y)=0\}$ is a singular planar algebraic curve. $ \qed $ \\

\begin{remark}
By Corollary 3.5, with the value $N=1$ we get the result that $\mathbb{C}[V,VU,VU^2+U]$ (which equals
to $I((X^{-1},X^2 Y-X))$) contains no counterexample to the Jacobian Conjecture. This was originally
proved by L. Makar-Limanov,(See \cite{rp,rp1}). He used in a clever way a grading technique on this
algebra giving the weights $\pm 1$ to the indeterminates $U$ and $V$ respectively. We recall that
Pinchuk's counterexample to the Real Jacobian Conjecture is contained in the real version $\mathbb{R}[V,VU,VU^2+U]$.
\end{remark}

\section{A necessary condition on the phantom curves for the surjectivity of the mapping}

\begin{theorem}
If $F\in\mathbb{C}[X,Y]^2$ satisfies $\det J_F(X,Y)\in\mathbb{C}^{\times}$, and
$\forall\,R\in R_0(F),\,\,\{ S_R(X,Y)=0\}\cap {\rm sing}(R)=\emptyset$, then $F(\mathbb{C}^2)=\mathbb{C}^2$.
\end{theorem}
\noindent
{\bf Proof.} \\
If $F\in {\rm Aut}(\mathbb{C}^2)$, then the claim is true. If $F\not\in {\rm Aut}(\mathbb{C}^2)$, then
$F$ is a counterexample to the Jacobian Conjecture. In this case it has a non empty geometric basis $R_0(F)$.
By pre composing $F$ with a suitable invertible linear mapping $L\,:\,\mathbb{C}^2\rightarrow\mathbb{C}^2$ we 
can achieve the situation that $F$ satisfies the $Y$-degree condition: $\deg F=\deg_Y P=\deg_Y Q$. This implies that
each $R\in R_0(F)$ could be chosen to have the following form: $R(X,Y)=(X^{-\alpha},X^{\beta}Y+X^{-\alpha}\phi(X))$,
$\alpha,\beta\in\mathbb{Z}^+$, $\alpha<\beta$, $\phi(X)\in\mathbb{C}[X]$, $\deg\phi<\alpha+\beta$ and the gcd
of the set of $X$-powers that effectively appear in $X^{\alpha+\beta}Y+\phi(X)$ equals to $1$.
Also if $H_R(X,Y)=0$ is an implicit representation of the $R$-irreducible component of $A(F)$, then by Theorem 3.1
$(H_R\circ G_R)(X,Y)=X^{\beta-\alpha}S_R(X,Y)$ where $S_R(X,Y)\in\mathbb{C}[X,Y]$. Finally by the assumptions we have
$\{ S_R(X,Y)=0\}\cap {\rm sing}(R)=\emptyset$. The $G_R$-pre-image of the $R$-irreducible component of the asymptotic 
variety $A(F)$ is the union of singular locus of $R$, ${\rm sing}(R)=\{ X=0\}$ and of the $R$-phantom curve 
$\{ S_R(X,Y)=0\}$. More accurately we have $ G_R^{-1}(\{ H_R(X,Y)=0\})=G_R^{-1}(G_R({\rm sing}(R)))={\rm sing}(R)\cup\{ S_R(X,Y)=0\}$.
The irreducible curves $\{ H_R(X,Y)=0\}$ are exactly the set of the asymptotic values of $F$, $A(F)$, when we take
all the rational mappings $R\in R_0(F)$. As is the tradition in complex analysis we call the asymptotic values of $F$ 
which do not belong to its image, the Picard exceptional values of $F$. We denote this set by ${\rm Picard}(F)$, and
$\forall\,R\in R_0(F)$ we denote the $R$-Picard exceptional values of $F$ by ${\rm Picard}_R(F)$. Thus: 
${\rm Picard}_R(F)={\rm Picard}(F)\cap \{ H_R(X,Y)=0\}$. Hence we have the representation: ${\rm Picard}(F)=
\bigcup_{R\in R_0(F)} {\rm Picard}_R(F)$. As is well known we have $\mathbb{C}^2-F(\mathbb{C}^2)={\rm Picard}(F)$, and
the finiteness $0\le |{\rm Picard}(F)|<\infty$. Our theorem is merely the assertion ${\rm Picard}(F)=\emptyset$, 
or, equivalently $\forall\,R\in R_0(F)$, ${\rm Picard}_R(F)=\emptyset$. We will prove this last assertion. Let us fix
an element in the geometric basis of $F$, $R\in R_0(F)$. By the above, the difference set $\{ H_R(X,Y)=0\}-
G_R(\{ S_R(X,Y)=0\})$ is a finite subset of the $R$-irreducible component $G_R(\{ X=0\})=\{ H_R(X,Y)=0\}$. If $F$
is not surjective, i.e. $F(\mathbb{C}^2)\ne\mathbb{C}^2$, then $\mathbb{C}^2-F(\mathbb{C}^2)$ is a finite set
which is composed exactly of these asymptotic values of $F$ which are the Picard exceptional values of $F$. We know that 
$\{ X=0\}\cap\{ S_R(X,Y)=0\}$ is empty. It follows that $R$ is defined on all the points of the $R$-phantom curve 
$\{ S_R(X,Y)=0\}$. Hence by the definition of the $R$-dual mapping of $F$ we have, $G_R(\{ S_R(X,Y)=0\})=F(R(\{ S_R(X,Y)=0\}))
\subseteq F(\mathbb{C}^2)$. Hence if $H_R(X,Y)=0$ contains Picard exceptional values of $F$, i.e. if ${\rm Picard}_R(F)\ne\emptyset$, 
then they form a subset of the finite set $\{ H_R(X,Y)=0\}-G_R(\{ S_R(X,Y)=0\})$. We conclude that if we prove that, 
$\{ H_R(X,Y)=0\}-G_R(\{ S_R(X,Y)=0\})\subseteq F(\mathbb{C}^2)$, then $\{ H_R(X,Y)=0\}$ contains no Picard exceptional values of $F$, i.e.
${\rm Picard}_R(F)=\emptyset$. Since $R$ is an arbitrary member of the geometric basis $R_0(F)$ of $F$ this
would imply that $F$ is surjective. Thus we now turn to prove that: $\{ H_R(X,Y)=0\}-G_R(\{ S_R(X,Y)=0\})
\subseteq F(\mathbb{C}^2)$. If $G_R$ is proper on an irreducible component $L$ of the $R$-phantom curve
$S_R(X,Y)=0$, then $G_R(L)=\{ H_R(X,Y)=0\}$ because $H_R(X,Y)=0$ is an irreducible component of $A(F)$.
Since we have $\{ X=0\}\cap L=\emptyset$, it follows that $R$ is defined on $L$ and so, $\{ H_R(X,Y)=0\}=
G_R(L)=F(R(L))\subseteq F(\mathbb{C}^2)$. If $G_R$ is not proper on any component of the $R$-phantom curve
$S_R(X,Y)=0$, then any such a component is an asymptotic tract of $G_R$. Let $\{(f(T),g(T))\,|\,T\in\mathbb{C}\}$
be a parametrization of the component $L$ of the $R$-phantom curve. Then, as explained above, $R$ is defined on $L$ and we
have, $R(L)=\{(f(T)^{-\alpha},f(T)^{\beta}g(T)+f(T)^{-\alpha}\phi(f(T)))\}$. There are two possibilities: (1) The curve
$\{(f(T)^{-\alpha},f(T)^{\beta}g(T)+f(T)^{-\alpha}\phi(f(T)))\}$ is an asymptotic tract of the mapping $F$. (2) This curve,
$\{(f(T)^{-\alpha},f(T)^{\beta}g(T)+f(T)^{-\alpha}\phi(f(T)))\}$ is bounded. \underline{Claim}: (1) is impossible. \\
\underline{A proof of the claim}: (1)$\Rightarrow\,f(T)\rightarrow 0$ and $g(T)$ stays bounded. This follows
because $F$ satisfies the $Y$-degree condition (see equation (*) in the proof of Theorem 2.1). But then the
component $L=\{(f(T),g(T))\}$ of the $R$-phantom curve is bounded. This is not possible $ \qed$ \\
It is worth giving \underline{a second proof of the claim}: As in the first proof $f(T)\rightarrow 0$, $g(T)$
stays bounded. On the other hand $(f(T),g(T))$ is a parametrization of $L$ which is a component of
$S_R(X,Y)=0$. Now we have the representation $S_R(X,Y)=e_R+X\cdot T_R(X,Y)$ for some $e_R\in\mathbb{C}^{\times}$
and some $T_R(X,Y)\in\mathbb{C}[X,Y]$ (because $\{S_R(X,Y)=0\}\cap\{ X=0\}=\emptyset$). Thus $S_R(f(T),g(T))\equiv 0$ which implies that $e_R+f(T)\cdot
T_R(f(T),g(T))\equiv 0$. But by $f(T)\rightarrow 0$ and $g(T)$ stays bounded we deduce that $e_R=0$ which
contradicts $e_R\in\mathbb{C}^{\times}$. $\qed $ \\
Thus only possibility (2) occurs. In this case $f(T)\rightarrow c\in\mathbb{C}^{\times}\cup\{\infty\}$. If
$c\in\mathbb{C}^{\times}$ then $g(T)\rightarrow c_1\in\mathbb{C}$ (otherwise the curve 
$R(L)=\{(f^{-\alpha},f^{\beta}g+f^{-\alpha}\phi(f))\}$
is not bounded). This implies that $L=\{(f(T),g(T))\}$ is a bounded curve. This is not possible. We deduce
that $f(T)\rightarrow\infty$ and $g(T)\rightarrow 0$ (otherwise the $Y$-coordinate of $R(L)$, $f^{\beta}g+f^{-\alpha}\phi(f)
\rightarrow\infty$ because $\deg \phi(X)<\alpha+\beta$). Thus $f(T)\rightarrow\infty$ and $g(T)\rightarrow 0$ in such
a manner that $f^{\beta}g+f^{-\alpha}\phi(f)\rightarrow d\in\mathbb{C}$. We deduce that the asymptotic value of $G_R$
along $L$ is $F(0,d)$ and in particular, it belongs to the image of $F$, $F(\mathbb{C}^2)$. \\
\underline{Conclusion}: The $R$-dual mapping, $G_R$ of $F$ maps the $R$-phantom curve $S_R(X,Y)=0$ into $F(\mathbb{C}^2)$
and moreover the asymptotic values of $G_R$ along the components of $S_R(X,Y)=0$ also belong to $F(\mathbb{C}^2)$.
In fact they belong to $\{ F(0,Y)\,|\,Y\in\mathbb{C}\}$. This proves that $\{ H_R(X,Y)=0\}-G_R(\{ S_R(X,Y)=0\})
\subseteq F(\mathbb{C}^2)$, and concludes the proof of the surjectivity of the Keller mapping $F$. $\qed $ \\

\section{More on the structure of the $R$-phantom curve and a type of a Picard's (small) Theorem}

Let $F(U,V)$ be a Keller mapping which is not a $\mathbb{C}^2$-automorphism. Then $R_0(F)\ne\emptyset$.
Let $R(X,Y)=(X^{-\alpha},X^{\beta}Y+X^{-\alpha}\Phi(X))\in R_0(F)$. We recall that $\alpha,\,\beta\in\mathbb{Z}^+$,
$\beta-\alpha-1>0$, the gcd of all the $X$-powers that effectively appear in $X^{\alpha+\beta}Y+\Phi(X)$
equals to $1$, where $\Phi(X)\in\mathbb{C}[X]$, $\deg\Phi<\alpha+\beta$. We can further assume that
$X^{-\alpha}\Phi(X)\in\mathbb{C}[X]$. If $H_R(X,Y)=0$ is the $R$-irreducible component of $A(F)$, then
$(H_R\circ F)\circ R(X,Y)=X^{\beta-\alpha}S_R(X,Y)$ (Theorem 3.1), where $S_R(X,Y)\in\mathbb{C}[X,Y]$
and $S_R(X,Y)=0$ is the $R$-phantom curve. We have $G_R^{-1}(\{ H_R(X,Y)=0\})=\{X=0\}\cup\{S_R(X,Y)=0\}$,
where $G_R=F\circ R$ (off $\{ X=0\}$) is the $R$-dual of $F$. Thus we obtain by differentiations:
$$
\left(\frac{\partial}{\partial V}(H_R\circ F)\right)\circ R(X,Y)=X^{-\alpha}\frac{\partial S_R}{\partial Y},
$$
$$
\left(\frac{\partial}{\partial U}(H_R\circ F)\right)\circ R(X,Y)=-\frac{X^{-\alpha}}{\alpha}
\left\{(\beta-\alpha)X^{\alpha+\beta}S_R(X,Y)+\right.
$$
$$
\left. +X^{\alpha+\beta+1}\cdot\frac{\partial S_R}{\partial X}-
(\beta X^{\alpha+\beta}Y-\alpha\Phi(X)+X\Phi'(X))\cdot\frac{\partial S_R}{\partial Y}\right\}.
$$
We are interested in the intersection points of the $R$-phantom curve and of ${\rm sing}(R)=\{ X=0\}$. Let
$(0,Y_0)$ be such a point. Then $S_R(0,Y_0)=0$. We can represent $S_R(X,Y)$ as follows: $S_R(X,Y)=f(Y)+X\cdot g(X,Y)$,
where $f(Y)\in\mathbb{C}[Y]$ and $g(X,Y)\in\mathbb{C}[X,Y]$. Then $S_R(0,Y_0)=0\Rightarrow f(Y_0)=0$.
By Corollary 3.2 (or Corollary 3.3) $G_R(0,Y_0)$ is a singularity of the $R$-irreducible component of $A(F)$:
$$
\frac{\partial H_R}{\partial U}(G_R(0,Y_0))=\frac{\partial H_R}{\partial V}(G_R(0,Y_0))=0.
$$
$$
\Rightarrow 0=\frac{\partial H_R}{\partial V}(G_R(0,Y_0))=\lim_{X\rightarrow 0}X^{-\alpha}\frac{\partial S_R}{\partial Y}(X,Y_0)=
\lim_{X\rightarrow 0}X^{-\alpha}(f'(Y_0)+X\cdot\frac{\partial g}{\partial Y}(X,Y_0))
$$
$$
\Rightarrow f'(Y_0)=0\wedge 0=\lim_{X\rightarrow 0}X^{-\alpha+1}\frac{\partial g}{\partial Y}(X,Y_0)
$$
$$
\Rightarrow f'(Y_0)=0\wedge \frac{\partial g}{\partial Y}(X,Y_0)=X^{\alpha}h(X,Y_0)
$$
$$
\Rightarrow f'(Y_0)=0\wedge g(X,Y)=h_3(X)+(Y-Y_0)^2h_2(X,Y)+X^{\alpha}h_1(X,Y),
$$
where $\deg_X(h_3+(Y-Y_0)^2h_2)<\alpha$
$$
\Rightarrow S_R(X,Y)=(Y-Y_0)^2[f_1(Y)+Xh_2(X,Y)]+X[H_3(X)+X^{\alpha}h_1(X,Y)],
$$
where $\deg_X h_2(X,Y),\,\deg h_3(X)<\alpha$. We denote $\Psi(X,Y)=f_1(Y)+Xh_2(X,Y)$, then
$$
S_R(X,Y)=(Y-Y_0)^2\Psi(X,Y)+X[h_3(X)+X^{\alpha}h_1(X,Y)],                           \eqno(*)
$$
where $X\not|\Psi(X,Y),\,\,\deg_X\Psi(X,Y)\le\alpha,\,\,\deg_X h_3(X)<\alpha$. By computing the derivatives of
$S_R(X,Y)$ (in (*)) and substituting $(X,Y)=(0,Y_0)$ we obtain
$$
\frac{\partial S_R}{\partial Y}(0,Y_0)=0,\,\,\,\frac{\partial S_R}{\partial X}(0,Y_0)=h_3(0).
$$
So far we have used the equation $(\partial H_R/\partial V)(G_R(0,Y_0))=0$. We now turn to the second component of
the gradient of $H_R$ at the singular point $G_R(0,Y_0)$.
$$
0=\frac{\partial H_R}{\partial U}(G_R(0,Y_0))=\lim_{X\rightarrow 0,\,\,Y\rightarrow Y_0}\left(-\frac{X^{-\alpha}}{\alpha}
\left\{(\beta-\alpha)X^{\alpha+\beta}S_R(X,Y)+\right.\right.
$$
$$
+\left.\left. X^{\alpha+\beta+1}\cdot\frac{\partial S_R}{\partial X}-(\beta X^{\alpha+\beta}Y-\alpha\Phi(X)+X\Phi'(X))
\cdot\frac{\partial S_R}{\partial Y}\right\}\right)=
$$
$$
=\left(-\frac{1}{\alpha}\right)\lim_{X\rightarrow 0,\,\,Y\rightarrow Y_0} X^{-\alpha}\left\{-(-\alpha\Phi(X)+X\Phi'(X))
[2(Y-Y_0)\Psi(X,Y)+(Y-Y_0)^2\frac{\partial \Psi}{\partial Y}]\right\}.
$$
We conclude that
$$
X^{-\alpha}(-\alpha\Phi(X)+X\Phi'(X))[2(Y-Y_0)\Psi(X,Y)+(Y-Y_0)^2\frac{\partial \Psi}{\partial Y}]\in\mathbb{C}[X,Y].
$$
Also by $X^{-\alpha}\Phi(X)\in\mathbb{C}[X]$ we have $X^{-\alpha+1}\Phi'(X)\in\mathbb{C}[X]$ and hence
$X^{-\alpha}(-\alpha\Phi(X)+X\Phi'(X))\in\mathbb{C}[X]$. Moreover $(X^{-\alpha}\Phi(X))'=X^{-\alpha-1}(-\alpha\Phi(X)+X\Phi'(X))\in\mathbb{C}[X]$
so $\lim_{X\rightarrow 0} X^{-\alpha}(-\alpha\Phi(X)+X\Phi'(X))=0$.

Let us suppose that the $R$-phantom curve $S_R(X,Y)=0$ intersects ${\rm sing}(R)$ in the set of points 
$(0,Y_j)$, $0\le j\le N-1$. Then by the above calculation we have for each $0\le j\le N-1$:
$$
S_R(X,Y)=(Y-Y_j)^2\Psi_j(X,Y)+X[h_{3j}(X)+X^{\alpha}h_{1j}(X,Y)]. 
$$
We substitute $X=0$ and obtain:
$$
(Y-Y_0)^2\Psi_0(0,Y)=(Y-Y_1)^2\Psi_1(0,Y)=\ldots=(Y-Y_{N-1})^2\Psi_{N-1}(0,Y).
$$
Since $Y_0,\ldots,Y_{N-1}$ are the total set of zeros of $S_R(0,Y)$ we obtain:

\begin{proposition}
$$
S_R(X,Y)=\Psi(X)\prod_{j=0}^{N-1}(Y-Y_j)^{2+\epsilon_j}+X[h_3(X)+X^{\alpha}h_1(X,Y)],        \eqno(**)
$$
$\Psi(0)\ne 0$, $\deg\Psi(X)\le\alpha$, $\epsilon_j\in\mathbb{Z}^+\cup\{ 0\}$, $\deg h_3(X)<\alpha$.
\end{proposition}

\begin{corollary}
Either all the intersection points ${\rm sing}(R)\cap\{ S_R(X,Y)=0\}$ are critical points of
$S_R(X,Y)$ or none is such a critical point.
\end{corollary}
\noindent
{\bf Proof.} \\
This follows by 
$$
\forall\,0\le j\le N-1,\,\,\,\frac{\partial S_R}{\partial Y}(0,Y_j)=0,\,\,\,\frac{\partial S_R}{\partial X}(0,Y_j)=h_3(0).
$$
$\qed $ \\
By Corollary 3.3 a point $(X_0,Y_0)\not\in {\rm sing}(R)$ is a singular point of $S_R(X,Y)=0$ iff 
$G_R(X_0,Y_0)$ is a singular point of $H_R(X,Y)=0$. If we substitute equation (**) of Proposition 5.1 into the
basic relation $(H_R\circ F)\circ R(X,Y)=X^{\beta-\alpha}S_R(X,Y)$ take $X\rightarrow 0$ and remember that
$\Psi(0)\ne 0$ we obtain the following estimate for $X\rightarrow 0$ and $Y$ fixed:
$$
(H_R\circ F)(X^{-\alpha},X^{\beta}Y+X^{-\alpha}\Phi(X))=\left\{\begin{array}{lll}
\Omega(X^{\beta-\alpha}) & {\rm for} & Y\not\in\{Y_0,\ldots,Y_{N-1}\} \\
O(X^{\beta-\alpha+1}) & {\rm for} & Y\in\{Y_0,\ldots,Y_{N-1}\}\end{array}\right.
$$
The $R$-phantom curve does not intersect ${\rm sing}(R)$ iff $\{Y_0,\ldots,Y_{N-1}\}=\emptyset$ and this is equivalent to:
$$
\lim_{X\rightarrow 0}\frac{(H_R\circ F)(X^{-\alpha},X^{\beta}Y+X^{-\alpha}\Phi(X))}{X^{\beta-\alpha}}=c\in\mathbb{C}^{\times}.
$$
This is equivalent to:
$$
\left\{\frac{\partial^{\beta-\alpha+1}}{\partial Y\partial X^{\beta-\alpha}}
((H_R\circ F)(X^{-\alpha},X^{\beta}Y+X^{-\alpha}\Phi(X)))\right\}_{X=0}=0.
$$
This is equivalent to:
$$
\frac{\partial^{\beta-\alpha}}{\partial X^{\beta-\alpha}}\left\{\left(\frac{\partial(H_R\circ F)}{\partial V}
(X^{-\alpha},X^{\beta}Y+X^{-\alpha}\Phi(X))\right)\cdot X^{\beta}\right\}=X\cdot h_5(X,Y),
$$
for some $h_5(X,Y)\in\mathbb{C}[X,Y]$. Here we think of $(H_R\circ F)(U,V)$. Since $H_R\circ F\in I(R)$ it follows that
$$
X^{\beta}\cdot\frac{\partial(H_R\circ F)}{\partial V}(X^{-\alpha},X^{\beta}Y+X^{-\alpha}\Phi(X))\in\mathbb{C}[X,Y],
$$
which is equivalent to
$$
X^{\alpha-1}\cdot\frac{\partial(H_R\circ F)}{\partial V}(X^{-\alpha},X^{\beta}Y+X^{-\alpha}\Phi(X))=h_6(X,Y)\in\mathbb{C}[X,Y].
$$
We proved the following:
\begin{theorem}
Let $F$ be a Keller mapping which is not a $\mathbb{C}^2$-automorphism. Let
$R(X,Y)=(X^{-\alpha},X^{\beta}Y+X^{-\alpha}\Phi(X))\in R_0(F)$. Then $\{ S_R(X,Y)=0\}\cap {\rm sing}(R)=\emptyset$ iff
$$
\frac{\partial(H_R\circ F)}{\partial V}(R)=\frac{h_6(X,Y)}{X^{\alpha-1}},
$$
for some $h_6(X,Y)\in\mathbb{C}[X,Y]$.
\end{theorem}

\begin{remark}
Theorem 4.1 says that the disjointness $\{ S_R(X,Y)=0\}\cap {\rm sing}(R)=\emptyset$
condition that appears in Theorem 5.3 is sufficient for the surjectivity of the Keller mapping $F$, i.e. for 
$F(\mathbb{C}^2)=\mathbb{C}^2$. At this point it looks as if the condition given in Theorem 5.3 that is equivalent
to  $\{ S_R(X,Y)=0\}\cap {\rm sing}(R)=\emptyset$ is improbable to hold for all Keller mappings. The reason is that apriori
we only have:
$$
X^{\beta}\cdot\left(\frac{\partial(H_R\circ F)}{\partial V}\right)\circ R\in\mathbb{C}[X,Y].
$$
This follows immediately from $H_R\circ F\in I(R)$. However according to Theorem 5.3 we would like to have:
$$
X^{\alpha-1}\cdot\left(\frac{\partial(H_R\circ F)}{\partial V}\right)\circ R\in\mathbb{C}[X,Y],
$$
which is far away from what we have. However, we recall that in the beginning of this section we noticed that:
$$
X^{\alpha}\cdot\left(\frac{\partial(H_R\circ F)}{\partial V}\right)\circ R=\frac{\partial S_R}{\partial Y}\in\mathbb{C}[X,Y].
$$
This looks promising. We only have a difference of $1$ between $\alpha$ and $\alpha-1$.
\end{remark}

\begin{theorem}
If $F\in\mathbb{C}[X,Y]^2$ satisfies $\det J_F(X,Y)\in\mathbb{C}^{\times}$, then 
$$
\mathbb{C}^2-F(\mathbb{C}^2)\subseteq \bigcup_{R\in R_0(F)} G_R({\rm sing}(R)\cap\{ S_R(X,Y)=0\}).
$$
Also 
$$
\mathbb{C}^2-F(\mathbb{C}^2)\subseteq \bigcup_{R\in R_0(F)} {\rm sing}(\{ H_R(X,Y)=0\}).
$$
\end{theorem}
\noindent
{\bf Proof.} \\
Each $R\in R_0(F)$ contributes the $R$-Picard exceptional values of $F$ from among 
the $G_R$ images of the intersection points of the two two curves ${\rm sing}(R)$ and the $R$-phantom curve.
These give the entire set of the Picard exceptional values of $F$. This proves the first equation. The second follows by 
Corollary 3.3 which implies that:
$$
\forall\,R\in R_0(F),\,\,G_R({\rm sing}(R)\cap\{ S_R(X,Y)=0\})\subseteq {\rm sing}(\{ H_R(X,Y)=0\}).
$$
$\qed $ \\

How large can the set of the Picard exceptional values of $F$ be?

\begin{theorem}
If $F\in\mathbb{C}[X,Y]^2$ satisfies $\det J_F(X,Y)\in\mathbb{C}^{\times}$, then 
$$
|\mathbb{C}^2-F(\mathbb{C}^2)|\le (\deg F)^3+(\deg F)^2-(\deg F).
$$
\end{theorem}
\noindent
{\bf Proof.} \\
Our starting point will be the result in Theorem 5.5:
$$
\mathbb{C}^2-F(\mathbb{C}^2)\subseteq \bigcup_{R\in R_0(F)} G_R({\rm sing}(R)\cap\{ S_R(X,Y)=0\}).
$$
By the Bezout Theorem $|{\rm sing}(R)\cap\{ S_R(X,Y)=0\}|\le\deg S_R(X,Y)$. In fact we have by Proposition 5.1
$|{\rm sing}(R)\cap\{ S_R(X,Y)=0\}|=|\{ S_R(0,Y)=0\}|=\deg_Y S_R(0,Y)$. We know that
$G_R^{-1}(G_R({\rm sing}(R)))={\rm sing}(R)\cup\{ S_R(X,Y)=0\}$ and we are led to consider 
$$
\deg G_R^{-1}(G_R(0,Y))-\deg {\rm sing}(R)=\deg G_R^{-1}(G_R(0,Y))-1.
$$
So that
$$
|\mathbb{C}^2-F(\mathbb{C}^2)|\le\sum_{R\in R_0(F)}(\deg G_R^{-1}(G_R(0,Y))-1).
$$
It follows that
$$
|\mathbb{C}^2-F(\mathbb{C}^2)|\le\sum_{R\in R_0(F)}(\deg G_R -1).
$$
We need to estimate $\deg G_R$. Off ${\rm sing}(R)$ we have:
$$
G_R(X,Y)=(F\circ R)(X,Y)=F(X^{-\alpha},X^{\beta}Y+X^{-\alpha}\Phi(X)).
$$
Let us take a coordinate of $F$, $P(U,V)=\sum_{i+j\le N}a_{ij}U^iV^j$. Then on composition with $R$ we get
$$
(P\circ R)(X,Y)=\sum_{i+j\le N}a_{ij}(X^{-\alpha})^i(X^{\beta}Y+X^{-\alpha}\Phi(X))^j.
$$
The degree of a generic monomial is $(\beta+1)j-i\alpha$, $i+j\le N$. So we are looking at $\max\{\beta+1)j-i\alpha\,|\,
i+j=N\}=(\beta+1)\cdot N$. We arrive at the estimate $\deg G_R\le (\beta+1)\cdot\deg F$ and hence
$$
|\mathbb{C}^2-F(\mathbb{C}^2)|\le\sum_{R\in R_0(F)}(\deg F\cdot ((\beta+1)\deg F-1).
$$
From this it follows that
$$
|\mathbb{C}^2-F(\mathbb{C}^2)|\le\deg F\cdot ((\deg F+1)\cdot\deg F-1).
$$
$\qed $ \\
\\
Theorem 5.6 gives a cubic estimate (in terms of the degree of $F$) for the size of the
set of the Picard exceptional values of $F$. 

\begin{remark}
The independent interest of bringing the results of this section is that we get a type of
a Picard's (small) Theorem for polynomial \'etale mappings $K^2\rightarrow K^2$.
This result mostly, does not require the field $K$ to be algebraically closed. It does
assume the special form of the elements in the geometric basis of $F$, i.e.
$$
R(X,Y)=(X^{-\alpha},X^{\beta}\cdot Y+X^{-\alpha}\Phi(X)).
$$
\end{remark}

\noindent
{\it Ronen Peretz \\
Department of Mathematics \\ Ben Gurion University of the Negev \\
Beer-Sheva , 84105 \\ Israel \\ E-mail: ronenp@math.bgu.ac.il} \\ 
 
\end{document}